\documentclass[12pt]{article}

\usepackage{amsmath}
\usepackage{amsfonts}

\newcommand{\bb}{\mbox{\boldmath$b$}}
\newcommand{\be}{\mbox{\boldmath$e$}}
\newcommand{\bh}{\mbox{\boldmath$h$}}

\newcommand{\bnu}{\mbox{\boldmath$\nu$}}

\newcommand{\bq}{\mbox{\boldmath$q$}}
\newcommand{\br}{\mbox{\boldmath$r$}}
\newcommand{\bs}{\mbox{\boldmath$s$}}

\newcommand{\bv}{\mbox{\boldmath$v$}}
\newcommand{\bw}{\mbox{\boldmath$w$}}
\newcommand{\bx}{\mbox{\boldmath$x$}}
\newcommand{\by}{\mbox{\boldmath$y$}}
\newcommand{\bze}{{\bf 0}}

\newcommand{\id}{\mbox{\rm I}}
\newcommand{\nul}{{\cal N}}
\newcommand{\qed}{\hspace{5mm}$\Box$}
\newcommand{\ran}{{\cal R}}
\newcommand{\rank}{\mbox{rank}}
\newcommand{\real}{{\bf R}}
\newcommand{\rj}{{\bf R}^j}
\newcommand{\rjj}{{\bf R}^{j \times j}}

\newcommand{\rn}{{\bf R}^n}

\newcommand{\rnn}{{\bf R}^{n \times n}}

\newcommand{\spn}{\mbox{span}}
\newcommand{\trans}{{\mbox{\scriptsize T}}}
\newtheorem{theorem}{Theorem}
\newtheorem{lemma}{Lemma}

\date{}

\begin{document}

\title{GMRES on singular systems revisited}

\author{Ken Hayami\footnote{National Institute of Informatics, and The Graduate University for Advanced Studies (SOKENDAI), 2-1-2, Hitotsubashi, Chiyoda-ku, Tokyo 101-8430, Japan, e-mail: hayami@nii.ac.jp} and Kota Sugihara\footnote{National Institute of Informatics, 2-1-2, Hitotsubashi, Chiyoda-ku, Tokyo 101-8430, Japan, e-mail: sugihara@nii.ac.jp} }

\maketitle

\begin{abstract}
In [Hayami K, Sugihara M. Numer Linear Algebra Appl. 2011; 18:449--469], the authors analyzed the convergence behaviour of the Generalized Minimal Residual (GMRES) method for the least squares problem 
$ \min_{\bx \in \rn} {\| \bb - A \bx \|_2}^2$, where $ A \in \rnn$ may be singular and $ \bb \in \rn$, by decomposing the algorithm into the range $ \ran(A) $ and its orthogonal complement $ \ran(A)^\perp $ components. However, we found that the proof of the fact that GMRES gives a least squares solution if $ \ran (A) = \ran(A^\trans) $ was not complete. In this paper, we will give a complete proof.
\end{abstract}

\section{Introduction}\label{intro}
In Hayami, Sugihara\cite{HS}, we showed in Theorem 2.6 that the Generalized Minimal Residual (GMRES) method of Saad, Schultz\cite{SS} gives a least squares solution to the
least squares problem
\begin{equation}
 \min_{\bx \in \rn} { \| \bb - A \bx \|_2}^2
 \label{lstsq}
\end{equation}
where $A \in \rnn $ may be singular, for all $\bb \in \rn$ and initial solution $\bx_0 \in \rn$ if and only if $\ran(A)=\ran(A^\trans)$, where $\ran(A)$ is the range space of $A$. The theorem had been proved by Brown and Walker\cite{BW}, but we gave an alternative proof by decomposing the algorithm into the $\ran(A)$ component and $\ran(A)^\perp$ component, thus giving a geometric interpretation to the range symmetry condition: $ \ran (A) = \ran(A^\trans) $. However, we later realized that the proof is not so obvious as we stated. In this paper, we will give a complete proof.

We assume exact arithmetic, and the following notations will be used.\\ 
\hspace*{6mm}$ V^\perp$: orthogonal complement of subspace $ V $ of 
$ \rn $. \\ 
For $ X \in \rnn $,\\
\hspace*{6mm}$ \ran(X) $: the range space of $ X $, i.e., the subspace spanned 
by the column vectors of $ X, $\\
\hspace*{6mm}$ \nul (X) $: the null space of $ X $, i.e., the subspace of 
vectors $ \bv \in \rn $ such that $ X \bv = \bze, $

\section{Convergence analysis of GMRES on singular systems}
\subsection{GMRES}
The GMRES method of Saad, Schultz\cite{SS} applied to (\ref{lstsq}) is given as follows.
\\ \\
\underline{{\bf GMRES}}\\ \\
Choose $\bx_0$. \\ 
$ \br_0 = \bb - A \bx_0 $ \\
$ \bv_1 = \br_0 / ||\br_0||_2 $ \\
For $ j = 1, 2, \cdots $ until satisfied do \\
\hspace{5mm} $ h_{i,j}=(\bv_i,A \bv_j )\hspace{4mm}(i=1,2,\ldots,j) $ \\ 
\hspace{5mm}
$ {\displaystyle \hat{\bv}_{j+1} = A \bv_j - \sum_{i=1}^j h_{i,j} \bv_i } $ \\
\hspace{5mm} $ h_{j+1,j} =||\hat{\bv}_{j+1} ||_2 $.
\hspace{4mm}If $ h_{j+1,j} =0, $ goto $ \ast $. \\
\hspace{5mm} $ \bv_{j+1} = \hat{\bv}_{j+1} / h_{j+1,j} $ \\
End do \\
$ \ast \, k:=j $ \\
Form the approximate solution \\ 
\hspace{5mm} $ \bx_k = \bx_0 + [ \bv_1,\ldots,\bv_k] \by_k $ \\
where $ \by = \by_k $ minimizes
$ ||\br_k||_2 = ||\beta\be_1-\overline{H}_k \by ||_2 $.\\

Here, $ \overline{H}_k = [ h_{i,j} ] \in {\bf R}^{ (k+1) \times k}$ is a Hessenberg matrix, i.e., $h_{i,j}=0$ for $i>j+1$. \, \,
$\beta={||\br_0 ||_2} $ \, and \, $ \be_1 = [1,0,\ldots,0]^\trans \in {\bf R}^{k+1} $. 
The method minimizes the residual norm $\|\br_k\|_2$, over the search space 
$ \bx_k = \bx_0 + \spn \{ \bv_1,\ldots,\bv_k \}$, where 
$\spn \{ \bv_1,\ldots,\bv_k \}
=\spn \{ \br_0, A\br_0,\ldots,A^{k-1} \br_0 \},$ and 
$ (\bv_i,\bv_j)=0\hspace{3mm}(i\neq j)$.
Let $V_j = \left[ v_1, \ldots, v_j \right]$. Then,
\begin{equation}
 AV_j =V_{j+1} \overline{H}_j 
 \label{AVO}
\end{equation}
holds.

The GMRES is said to break down when $ h_{j+1,j} = 0$. Then,
\begin{equation}
  A V_j=V_j H_j 
\label{AV}
\end{equation}
 holds, where $H_j \in {\bf R}^{ j \times j}$ consists of the firet $j$ rows of $\overline{H}_j$.

When $ A$ is nonsingular, the iterates of GMRES converges to the 
solution for all $\bb, \bx_0 \in \rn $ within at most $ n $ steps in exact 
arithmetic \cite{SS}.

For the general case when $A$ may be singular, we define the following.

\subsection{A geometrical framework}
In this section we will begin by giving geometric interpretations to the 
conditions $\nul (A)=\nul (A^\trans)$ and 
$\ran (A) \cap \nul (A)=\{\bze \}$. 
This is done by decomposing the space $\rn$ into $\ran (A)$ and 
${\ran (A)}^\perp $.

Let $ \rank A = \dim \ran (A) = r > 0, $ and
\begin{equation}
\hspace*{-1.5cm}
\bq_1,\ldots,\bq_r: \mbox{\rm orthonormal basis of } \ran (A),
 \label{q1r}
\end{equation}
\begin{equation}
\hspace*{-1cm}
 \bq_{r+1},\ldots,\bq_n: \mbox{\rm orthonormal basis of } {\ran (A)}^\perp,
\end{equation}
\begin{equation}
\hspace*{-3.7cm}
 Q_1:= [ \, \bq_1, \ldots, \bq_r ] \in {\bf R}^{n \times r},
\end{equation}
\begin{equation}
\hspace*{-2.7cm}
 Q_2:= [ \, \bq_{r+1}, \ldots, \bq_n ] \in {\bf R}^{n \times (n-r)},
\end{equation}
so that,
\begin{equation}
\hspace*{-4.3cm}
 Q:= [ \, Q_1, Q_2 ] \in \rnn
 \label{Q}
\end{equation}
is an orthogonal matrix satisfying
\begin{equation}
\hspace*{-1mm}
 Q^\trans Q = Q Q^\trans = \id_n ,
 \label{qq}
\end{equation}
where $ \id_n $ is the identity matrix of order $n$.

Orthogonal transformation of the coefficient matrix $A$ using $Q$ gives
\begin{equation}
   \tilde{A}:= Q^\trans A Q = \left[
                                 \begin{array}{cc}
                                    {Q_1}^\trans A Q_1 & {Q_1}^\trans A Q_2\\
                                                 0     &       0
                                 \end{array}
                              \right]      
             = \left[
                  \begin{array}{cc}
                     A_{11} & A_{12}\\
                      0     &  0   
                  \end{array}
               \right],
 \label{qtaq}
\end{equation}
since $ {Q_2}^\trans A Q = 0 $. 
Here, $A_{11}:= {Q_1}^\trans A Q_1 $ and $ A_{12}:= {Q_1}^\trans A Q_2 $. 

In Hayami, Sugihara\cite{HS} we derived the following properties concerning the sub-matrices 
$ A_{11} $ and $ A_{12}$ in (\ref{qtaq}).
\begin{theorem} \label{th:a11} 
$A_{11}: \mbox{nonsingular} \, \Longleftrightarrow \, 
\ran (A) \cap \nul (A) = \{ \bze \} .$ 
\end{theorem} 

\begin{lemma} \label{lem:a11r}$ A_{12} =0
\Longrightarrow A_{11}: \mbox{nonsingular}$
\end{lemma}

\begin{theorem} \label{th:a12}
$A_{12} = 0 \, \Longleftrightarrow \, 
\ran(A) = \ran (A^\trans) \, \Longleftrightarrow \, 
\nul(A) = \nul(A^\trans).$ 
\end{theorem}

Now we will consider decomposing iterative algorithms into the $\ran(A)$ and 
$\ran(A)^\perp$ components. In order to do so, we will use the transformation 
\[ 
   \tilde{\bv}:= Q^\trans \bv = [ Q_1 , Q_2 ]^\trans \bv
               = \left[ \begin{array}{c}
                          {Q_1}^\trans \bv \\
                          {Q_2}^\trans \bv    
                        \end{array}
                 \right]
               = \left[ \begin{array}{c}
                           \bv^1 \\
                           \bv^2    
                        \end{array}
                  \right], 
\]
\[
    \bv = Q \tilde{\bv} = [ Q_1 , Q_2 ] \left[ \begin{array}{c}
                                                 \bv^1 \\
                                                 \bv^2    
                                               \end{array}
                                        \right] = Q_1 \bv^1 + Q_2 \bv^2 ,
\]
cf. (\ref{q1r})-(\ref{qq}), to decompose a vector variable $\bv$ in the 
algorithm. Here, $ \bv^1 $ corresponds to the $ \ran (A) $ 
component $ Q_1 \bv^1 $ of $ \bv $, and $ \bv^2 $ corresponds to the 
$ \ran (A)^\perp $ component $ Q_2 \bv^2 $ of $ \bv $.

For instance, the residual vector $ \br:= \bb - A \bx $ is 
transformed into 
\[ \tilde{\br}:= Q^\trans \br = Q^\trans \bb - Q^\trans A Q (Q^\trans \bx), \]
or
\[  \left[ \begin{array}{c}
              \br^1 \\
              \br^2    
            \end{array}
     \right] 
     = \left[ \begin{array}{c}
                  \bb^1 \\
                  \bb^2    
              \end{array}
       \right]
       - \left[ \begin{array}{cc}
                   A_{11} & A_{12} \\
                    0     &  0   
                \end{array}
         \right]
         \left[ \begin{array}{c}
                    \bx^1 \\
                    \bx^2    
                 \end{array}
         \right],
\]
i.e.,
\begin{equation}
   \begin{array}{lll}
       \br^1 & = & \bb^1 - A_{11} \bx^1 - A_{12} \bx^2 \\ 
       \br^2 & = & \bb^2 .
   \end{array}
 \label{r1r2}
\end{equation}
Hence, in the least squares problem (\ref{lstsq}), we have
\begin{equation}
  { \| \bb - A \bx \|_2 }^2 = { \| \br \|_2 }^2 = { || \tilde{\br} ||_2 }^2
     = { \| \br^1 \|_2 }^2 + { \| \bb^2 \|_2 }^2 .
 \label{ineq}
\end{equation}

Note that it is not necessary to compute $Q$ or to decompose the algorithm into the $\ran(A)$ and $\ran(A)^\perp$ components in practice. It is only for the theoretical analysis.
\subsection{Decomposition of GMRES}
Based on the above geometric framework, we will analyze GMRES 
for the case when $ A$ is singular, by decomposing it into the $\ran(A)$ 
component and the $\ran(A)^\perp$ component as follows. \\ \\
\underline{\bf Decomposed GMRES (general case)} 
\begin{equation}
\begin{array}{lll}
\underline{ {\cal R}(A) \, \mbox{component} } & \hspace{5mm} & 
\underline{ {{\cal R}(A)}^\perp \, \mbox{component} } \\ \\
\bb^1 = {Q_1}^\trans \bb & \hspace{5mm} & \bb^2 = {Q_2}^\trans \bb \\ \\
\multicolumn{3}{l}{\mbox{Choose} \, \, \bx_0} \\ \\
\bx_0^1 = {Q_1}^\trans \bx_0 & \hspace{5mm} & 
\bx_0^2 = {Q_2}^\trans \bx_0 \\ \\
\br_0^1 = \bb^1 - A_{11} \bx_0^1 - A_{12} \bx_0^2 & 
\hspace{5mm} & \br_0^2 = \bb^2 \\
||\br_0||_2 = \sqrt{ {||\br_0^1||_2 }^2 + { ||\bb^2||_2 }^2 } \\ \\
\bv_1^1 = \br_0^1 /||\br_0||_2 & \hspace{5mm} & 
\bv_1^2 = \bb^2/||\br_0||_2 \\ \\
\multicolumn{3}{l}{ \mbox{For} \, \, j = 1, 2, \ldots \,  \mbox{until satisfied do} } \\ \\
\multicolumn{3}{l}
{\hspace{1cm} h_{i,j}=(\bv_i^1, A_{11} \bv_j^1 + A_{12} \bv_j^2 )
\hspace{4mm}(i=1,2,\ldots,j)} \nonumber \\ \\
\hspace{1cm}{\displaystyle \hat{\bv}_{j+1}^1 
  = A_{11} \bv_j^1 + A_{12} \bv_j^2 
     - \sum_{i=1}^j h_{i,j} \bv_i^1 }& \hspace{5mm} & 
     {\displaystyle \hat{\bv}_{j+1}^2 = - \sum_{i=1}^j h_{i,j} \bv_i^2 } \\ \\
\multicolumn{3}{l}
{\hspace{1cm}
h_{j+1,j} 
 = \sqrt{ { ||\hat{\bv}_{j+1}^1 ||_2 }^2 + { ||\hat{\bv}_{j+1}^2 ||_2 }^2 }.
\hspace{5mm}\mbox{If} \, \, h_{j+1,j} =0, \, \, \mbox{goto} \, \ast. }\\ \\
\hspace{1cm}\bv_{j+1}^1 = \hat{\bv}_{j+1}^1 /h_{j+1,j} & \hspace{5mm} &
\bv_{j+1}^2 = \hat{\bv}_{j+1}^2 /h_{j+1,j} \\ \\
\mbox{End do} & & \\ \\
\ast \, k:=j & & 
\end{array}
\end{equation}
\clearpage
\begin{equation}
\begin{array}{lll}
\multicolumn{3}{l}{\mbox{Form the approximate solution}} \\ \\
\hspace{5mm}\bx^1_k=\bx_0^1 + [ \bv_1^1,\ldots,\bv_k^1] \, \by_k& \hspace{5mm}&
\bx^2_k=\bx_0^2 + [ \bv_1^2,\ldots,\bv_k^2] \, \by_k \\ \\
\multicolumn{3}{l}
{\mbox{where} \, \, \by=\by_k \, \, \mbox{minimizes} \, \,
||\br_k||_2 = ||\beta\be_1-\overline{H}_k \by ||_2 .}
\end{array}
\label{decgmres}
\end{equation}

From the above decomposed form of GMRES, we obtain
\begin{equation}
         \left[ \begin{array}{cc}
                   A_{11} & A_{12} \\
                    0     &  0   
                \end{array}
         \right]
   \left[ \begin{array}{c}
              V^1_j \\
              V^2_j
           \end{array}
        \right]
=  \left[ \begin{array}{c}
              V^1_{j+1} \\
              V^2_{j+1}
           \end{array}
        \right]
   \overline{H}_j ,
 \label{AVOD}
\end{equation}
which is equivalent to (\ref{AVO}),
where $\left[ V^l_j \right] =\left[ v^l_1, \ldots, v^l_j \right] \: (l=1,2)$.

When $h_{j+1,j}=0$, (\ref{AVOD}) becomes
\[  \left[ \begin{array}{cc}
                   A_{11} & A_{12} \\
                    0     &  0   
                \end{array}
         \right]
   \left[ \begin{array}{c}
              V^1_j \\
              V^2_j
           \end{array}
        \right]
=  \left[ \begin{array}{c}
              V^1_j \\
              V^2_j
           \end{array}
        \right]
   H_j ,
\]
which is equivalent to (\ref{AV}).

In passing, when the system is consistent, i.e. $\bb \in \ran(A)$, then $\bb^2=Q_2^\trans \bb =\bze$. Hence, in the $\ran(A)^\perp$ component of the above decomposed algorithm, 
$\br_0^2=\bb^2=\bze, \: \bv_1^2=\bze$. Thus, $\hat{\bv}_{l}^2=\bze$ and $\bv_l^2=\bze$ for $l=1,\ldots,j+1$.  Hence, $V_j^2=0, \: V_{j+1}^2=0$. Thus, (\ref{AVOD}) reduces to
\[ A_{11} V_j^1 = V_{j+1}^1 \overline{H}_j . \]
(See section 2.5 of Hayami, Sugihara\cite{HS}.)

Returning to the general case when the system may be inconsistent, in Theorem \ref{th:a12} we gave a geometric interpretation: $A_{12}=0$ to the
condition: $\nul(A)=\nul(A^\trans)$. 
Now it is important to notice that if $A_{12}=0$ holds, the decomposed GMRES 
further simplifies as follows. \\ \\
\underline{\bf Decomposed GMRES (Case $\nul(A)=\nul (A^\trans)$)}
\begin{equation}
\begin{array}{lll}
\underline{ {\cal R}(A) \, \mbox{component} } & \hspace{5mm} & 
\underline{ {{\cal R}(A)}^\perp \, \mbox{component} } \\ \\
\bb^1 = {Q_1}^\trans \bb & \hspace{5mm} & \bb^2 = {Q_2}^\trans \bb \\ \\
\multicolumn{3}{l}{\mbox{Choose} \, \, \bx_0} \\ \\
\bx_0^1 = {Q_1}^\trans \bx_0 & \hspace{5mm} & 
\bx_0^2 = {Q_2}^\trans \bx_0 \\ \\
\br_0^1 = \bb^1 - A_{11} \bx_0^1 & 
\hspace{5mm} & \br_0^2 = \bb^2 \\ \\
||\br_0||_2 = \sqrt{ {||\br_0^1||_2 }^2 + { ||\bb^2||_2 }^2 } \\ \\
\bv_1^1 = \br_0^1 /||\br_0||_2 & \hspace{5mm} & 
\bv_1^2 = \bb^2/||\br_0||_2 \\ \\
\multicolumn{3}{l}{ \mbox{For} \, \, j = 1, 2, \ldots \,  \mbox{until satisfied do} } \\ \\
\multicolumn{3}{l}
{\hspace{1cm} h_{i,j}=(\bv_i^1, A_{11} \bv_j^1 )
\hspace{4mm}(i=1,2,\ldots,j)} \nonumber \\ \\
\hspace{1cm}{\displaystyle \hat{\bv}_{j+1}^1 
  = A_{11} \bv_j^1 - \sum_{i=1}^j h_{i,j} \bv_i^1 }& \hspace{5mm} & 
     {\displaystyle \hat{\bv}_{j+1}^2 = - \sum_{i=1}^j h_{i,j} \bv_i^2 } \\ \\
\multicolumn{3}{l}
{\hspace{1cm}
h_{j+1,j} 
 = \sqrt{ { ||\hat{\bv}_{j+1}^1 ||_2 }^2 + { ||\hat{\bv}_{j+1}^2 ||_2 }^2 }.
\hspace{5mm}\mbox{If} \, \, h_{j+1,j} =0, \, \, \mbox{goto} \, \ast. } \\ \\
\hspace{1cm}\bv_{j+1}^1 = \hat{\bv}_{j+1}^1 /h_{j+1,j} & \hspace{5mm} &
\bv_{j+1}^2 = \hat{\bv}_{j+1}^2 /h_{j+1,j} \\ \\
\mbox{End do} & & \\ \\
\ast \, k:=j & & \\ \\
\multicolumn{3}{l}{\mbox{Form the approximate solution}}\\ \\
\hspace{5mm}\bx^1_k=\bx_0^1 + [ \bv_1^1,\ldots,\bv_k^1] \, \by_k& \hspace{5mm}&
\bx^2_k=\bx_0^2 + [ \bv_1^2,\ldots,\bv_k^2] \, \by_k\\ \\
\multicolumn{3}{l}
{\mbox{where} \, \, \by=\by_k \, \, \mbox{minimizes} \, \,
||\br_k||_2 = ||\beta\be_1-\overline{H}_k \by ||_2 .}
\end{array}
\label{decgmresa120}
\end{equation}

Then, (\ref{AVOD}) simplifies to
\begin{eqnarray}
  A_{11} V^1_j & = & V^1_{j+1} \overline{H}_j  \label{avvj1} \\
   0           & = & V^2_{j+1} \overline{H}_j . \nonumber
\end{eqnarray}
If further, $h_{j+1,j}=0$, we have
\begin{eqnarray}
  A_{ 11}V^1_j & = & V^1_{j} H_j  \label{avvj} \\
   0           & = & V^2_{j} H_j . \nonumber
\end{eqnarray}

Note here that the $\ran(A)$ component of GMRES is ``essentially equivalent'' to GMRES applied to 
$A_{11} \bx^1 = \bb^1$, except for the scaling factors for $\bv_j^1$. Note also that, from Lemma \ref{lem:a11r}, 
$ A_{12} = 0$ implies that $A_{11}$ is nonsingular. 
From these observations, we concluded in Hayami, Sugihara\cite{HS} (Section 2.3, p. 454) that if $ A_{12} = 0$, ``arguments similar to Saad, Schultz\cite{SS} for GMRES on nonsingular systems imply that GMRES gives a least-squares solution for all $\bb$ and $\bx_0$''. 

However, we later found that the proof is not so obvious. The difficulty is that, although the Krylov basis $V_1=\left[ \bv_1, \ldots, \bv_j \right]$ at step $j$ of the GMRES is orthonormal, the corresponding $\ran(A)$ component vecors $V^1_j=\left[ \bv^1_1, \ldots, \bv^1_j \right]$ are not necessarily orthogonal, and it is not even obvious that they are linearly independent.  In the following, we give a complete proof of the statement. See also Sugihara, Hayami, Zheng\cite{SHZ}, Theorem 1 for a related proof for the right-preconditioned MINRES method for symmetric singular systems.

First, we observe the following.
\begin{lemma}  \label{lem:b2}
In the GMRES method,
if $ \br_0 \ne \bze,  \; h_{i+1,i}\ne 0 \hspace{2mm} (1 \le i \le j-1)$, then $\bv_i^2=c_i \bb^2 \: (i=1,\ldots,j)$, i.e. all the $\ran(A)^\perp$ components $\bv_i^2 (i=1,\ldots,j)$ are parallel to $\bb^2$.
\end{lemma}
{\it Proof:} 
From the above Decomposed GMRES(general case) (\ref{decgmres}), \\ $\bv^2_1=\bb^2/||\br_0\|_2=c_1 \bb^2$. 
Since $\hat{\bv}_{j+1}^2 = - {\displaystyle \sum_{i=1}^j h_{i,j} \bv_i^2 }$ and $\bv_{j+1}^2 = \hat{\bv}_{j+1}^2 /h_{j+1,j} $, by induction, we have  $\bv_i^2=c_i \bb^2 \: (i=1,\ldots,j)$. \qed

Next, we prove the followi,ng.
\begin{theorem} \label{theorem:rv1j}
In the GMRES method,
assume $ \br_0 \ne \bze, \; h_{i+1,i}\ne 0 \; (1 \le i \le j-1)$ hold. If $\bb \in \ran(A) \: (\bb^2=0)$, then $\rank V_1^j=j$. If $\bb \notin \ran(A) \: (\bb^2 \ne \bze)$, then $\rank V^1_j=j-1$ or $j$.
\end{theorem}
{\it Proof:} 
When $\bb \in \ran(A) \: (\bb^2=0)$, from Lemma \ref{lem:b2}, 
\[ \tilde{V}_j=Q^T V_j = \left[ \begin{array}{c}
                      \bv_1^1, \ldots, \bv_j^1 \\
                      \bze, \ldots, \bze
                    \end{array}
           \right] .
\]
Hence, $\rank V^1_j = \rank V_j = j$. 

When $\bb \notin \ran(A) \: (\bb^2 \ne \bze)$, for $j=1$, $\rank V_1^1=\rank \left[ \bv_1^1 \right] = 0$ or $1$, depending on whether $\bv_1^1=\bze$ or $\bv_1^1 \ne \bze$. 

Let $j \ge 2$. From Lemma \ref{lem:b2}, and $ c_1 = 1/||\br_0 \| \ne 0$, we have
\[ \tilde{V}_j=Q^T V_j = \left[ \begin{array}{ccc}
                                                 \bv^1_1, & \ldots, & \bv^1_j  \\
                                                 c_1 \bb^2, & \ldots, & c_j \bb^2
                                              \end{array}
                                       \right]
 =  \left[ \begin{array}{cccc}
                      {\bv^1_1}',&{\bv^1_2}'& \ldots,&{\bv^1_j}' \\
                       \bb^2,&\bze,&\ldots,&\bze
             \end{array}
    \right] S^{-1} ,
\]
where
\[ S=  \left[ \begin{array}{cccc}
                     1/c_1 & -c_2/c_1 & \cdots & -c_j/c_1 \\
                              &  1         & \cdots & 0 \\
                              &             & \ddots & 0 \\
    $\mbox{\Large 0}$ &             &           & 1
                 \end{array}
        \right] \in \rjj 
\]
is nonsingular, and $ {\bv^1_i}'=\bv^1_i/c_1 \: (i=1,\ldots,j)$. Therefore, 
\[ \rank \left[ \begin{array}{cccc}
                      {\bv^1_1}',&{\bv^1_2}'& \ldots,&{\bv^1_j}' \\
                       \bb_2,&\bze,&\ldots,&\bze
             \end{array}
    \right] 
   = \rank V_j = j . \]
Then, $\rank \left[ {\bv_2}',\ldots,{\bv_j}' \right] =j-1$, since if $ \rank \left[ {\bv_2}',\ldots,{\bv_j}' \right] < j-1$, then 
\[ \rank  \left[ \begin{array}{cccc}
                      {\bv^1_1}',&{\bv^1_2}'& \ldots,&{\bv^1_j}' \\
                       \bb_2,&\bze,&\ldots,&\bze
             \end{array}
    \right] 
   < j . \]
Hence, $\rank \left[ \bv_1^1,\ldots \bv_j^1 \right] = \rank \left[ {\bv_1^1}',\ldots,{\bv^1_j}' \right] = j-1$ or $j$. 
\qed

Note that Lemma \ref{lem:b2} and Theorem \ref{theorem:rv1j} hold without assuming $A_{12}=0$.

Next, we prove the following, which corresponds to the sufficiency of the condition in Theorem 2.6 of Hayami, Sugihara\cite{HS}. 
\begin{theorem}
Assume $A_{12}=0$. Then, GMRES determines a least squares solution of (\ref{lstsq}) for all $\bb, \bx_0 \in \rn$.
\end{theorem}
{\it Proof:}
If $\br_0 = \bze$, a (least squares) solution to (\ref{lstsq}) is obtained. Assume $\br_0 \ne \bze$. 

Assume $ \bb \in \ran (A) $. Then, from Theorem \ref{theorem:rv1j}, $\rank V_1^j=j$. Since $\rank V_1^j \le r=\rank A$, there exists a $j\le r$, such that $ h_{i+1,i}\ne 0 \; (1 \le i \le j-1), \; h_{j+1,j} = 0$. Then from (\ref{avvj}), we have $A_{11} V_j^1 = V_j^1 H_j $. Since $A_{11}$ is nonsingular, $ \rank A_{11}V_1^j=j$. Then, $j=\rank V_j^1 H_j \le \min (j, \rank H_j )$, where $\rank H_j \le j$. Hence, $\rank H_j=j$, and $H_j$ is nonsingular. Note that
\begin{equation}
 \begin{array}{lll}
 \br_j^1 & = & \bb^1-A_{11} \bx_j^1 = \bb^1-A_{11} \left( \bx_0^1+V_j^1 \by_j \right) =  \br_0^1 - A_{11}V_j^1 \by_j  \\
           & = & \beta \bv_1^1 - V_j^1 H_j \by_j = V_j^1 \left( \beta \be_1 - H_j \by_j \right),  
 \end{array}
  \label{rj}
\end{equation}
where $\be_1 =(1,0,\ldots,0)^\trans \in \rj$. 
Hence, a least squares solution is obtained at step $j \; (j \le r) $ for $\by_j=\beta {H_j}^{-1} \be_1$, for which $\br_j^1=\bze.$

Next, assume $\bb \notin \ran(A)$. Then, in the proof of Theorem \ref{theorem:rv1j}, 
$ \rank A = r \ge \rank V_j^1 = j $ or $j-1$, which implies that there exists $j \le r+1$ such that $ h_{i+1,i}\ne 0 \hspace{2mm} (1 \le i \le j-1), \; h_{j+1,j}=0$.

(As in Point a and b in the proof of Theorem 1 in Sugihara et al.\cite{SHZ}), since $ V^2_{j} H_j = 0 $ from (\ref{avvj}), 
if $ H_j $ is nonsingular, $ V^2_j = [ \bv_1^2, \ldots, \bv_j^2 ] = 0$.
However, since $ \bb \notin \ran (A), \; \bb^2 \ne \bze$,
so that $ \bv_1^2 = \bb^2 / \| \br_0 \|_2 \ne \bze$.
Hence, $ H_j $ is singular, and there exists $\bw \ne \bze $ such that $ H_j \bw = \bze$.
(In fact, $\rank H_j = j-1$, since $h_{i+1,i} \ne 0 \; (1 \le i \le j-1)$.)
Then, from (\ref{avvj}), $ V_j^1 H_j \bw = A_{11} V_j^1 \bw = \bze$.
Since $ A_{11} $ is nosingular, $ V_j^1 \bw = \bze, \; \bw \ne \bze$. 
Hence, $ \rank V_j^1 = j-1$.
Then, a least squares solution is obtained at step $j$ if and only if $ H_j \by_j - \beta \be_1 \in \nul (V_j^1)$.
Since $\rank V_j^1 + \dim \nul (V_j^1)=j, \; \dim \nul(V_j^1)=1$. Let $ \nul(V_1^j) = \left\{ c \, \bnu^j \right\}$, where $c\in\real, \; \bnu \ne \bze \in \rj$. Let 
\[ \bnu = \left[ \begin{array}{c}
                          \nu_1 \\
                          \bnu_2
                       \end{array}
               \right] \ne \bze \in \real^j, \; \nu_1 \in \real, \; \bnu_2 \in {\bf R}^{j-1}, \; \mbox{and} \;
H_j = \left[ \begin{array}{cc}
                      {\bh_{11}}^\trans & h_{1j} \\
                      H_{21} & \bh_{22}
                    \end{array}
            \right],
\]
where $ {\bh_{11} }^\trans = \left[ h_{11},\ldots, h_{1,j-1} \right]$, 
\[ \begin{array}{ccc}
H_{21} = \left[ \begin{array}{ccc}
                            h_{21} & \cdots & h_{2,j-1} \\
                                     & \ddots & \vdots \\
                             $\mbox{\Large 0}$       &           & h_{j,j-1}
                          \end{array}
                  \right] & \mbox{and} & \bh_{22}=\left[ \begin{array}{c}
                                                                         h_{2j} \\
                                                                         \vdots \\
                                                                         h_{jj}
                                                                       \end{array}
                                                              \right].
    \end{array}
\]
where $ H_{21} $ is nonsingular since
\begin{equation}
    h_{i+1,i}\ne 0 \hspace{3mm} (1 \le i \le j-1).
 \label{hj,j-1}
\end{equation}
Note the following:
\begin{equation*}
 \begin{array}{ll}
 & \mbox{A least squares solution is obtained at step $j$ } \\
 \Longleftrightarrow & \exists \by \mbox{ such that } H_j \by - \beta \be_1 = c \bnu \\
 \Longleftrightarrow & \exists \by_1, y_j  \mbox{ such that } 
                                                              \left\{ \begin{array}{ccc}
                                                                           {\bh_{11} }^\trans \by_1 + h_{1j} y_j & = & \beta + c \nu_1 \\
                                                                           H_{21} \by_1 + y_j \bh_{22} & = & c \bnu_2 
                                                                      \end{array} 
                                                              \right. \\
 \Longleftrightarrow  & \left( h_{1j} - {\bh_{11} }^\trans { H_{21} }^{-1} \bh_{22} \right) y_j 
                                  = \beta + c \left( \nu_1 - {\bh_{11} }^\trans { H_{21} }^{-1} \bnu_2 \right), 
 \end{array}
\end{equation*}
where 
\begin{equation*}
 \by = \begin{array}{ccc}
            \left[ \begin{array}{c}
                        \by_1 \\
                        y_j
                    \end{array}
            \right]
& \mbox{and}
& \by_1 = \left[ \begin{array}{c}
                       y_1 \\
                      \vdots \\
                       y_{j-1}
               \end{array}
  \right] \in {\bf R}^{j-1}.
        \end{array}
\end{equation*}
Here note that
\begin{equation*}
 \left[ \begin{array}{cc}
            {\rm I} & \bze \\
            - {\bh_{11} }^\trans & 1
         \end{array}
 \right]
 \left[ \begin{array}{cc}
            {H_{21} }^{-1} & \bze \\
            {\bze }^\trans & 1
         \end{array}
\right]
 \left[ \begin{array}{cc}
            \bze & {\rm I} \\
             1 & {\bze }^\trans
         \end{array}
\right]
 \left[ \begin{array}{cc}
             {\bh_{11} }^\trans & h_{1j} \\
             H_{21} & \bh_{22}
         \end{array}
\right] =
 \left[ \begin{array}{cl}
             {\rm I } & {H_{21} }^{-1} \bh_{22} \\
             {\bze }^\trans & h_{1j} -  {\bh_{11} }^\trans { H_{21} }^{-1} \bh_{22}
        \end{array}
 \right].   
\end{equation*}
Since $ \det H_j = 0$, \: $ h_{1j} -  {\bh_{11} }^\trans { H_{21} }^{-1} \bh_{22} = 0$. Thus,
\begin{equation*}
 \begin{array}{ll}
 & \mbox{A least squares solution is obtained at step $j$ } \\
 \Longleftrightarrow & \beta = c \left( \nu_1 - {\bh_{11} }^\trans { H_{21} }^{-1} \bnu_2 \right) \\
 \Longleftrightarrow & \nu_1 - {\bh_{11} }^\trans { H_{21} }^{-1} \bnu_2 \ne 0
 \end{array}
\end{equation*}
since $ \beta \ne 0$. Hence, if $ \nu_1 - {\bh_{11} }^\trans { H_{21} }^{-1} \bnu_2 \ne 0 $, a least squares solution is obtained at step $j$.
If $ \nu_1 - {\bh_{11} }^\trans { H_{21} }^{-1} \bnu_2 = 0$, a least squares solution is not obtained at step $j$.
Note that
\begin{equation*}
 \left[ \begin{array}{cc}
            {\rm I} & \bze \\
            - {\bh_{11} }^\trans & 1
         \end{array}
 \right]
 \left[ \begin{array}{cc}
            {H_{21} }^{-1} & \bze \\
            {\bze }^\trans & 1
         \end{array}
\right]
 \left[ \begin{array}{cc}
            \bze & {\rm I} \\
             1 & {\bze }^\trans
         \end{array}
\right]
 \left[ \begin{array}{cc}
             {\bh_{11} }^\trans & \nu_1 \\
             H_{21} & \bnu_2
         \end{array}
\right] =
 \left[ \begin{array}{cl}
             {\rm I } & {H_{21} }^{-1} \bnu_2 \\
             {\bze }^\trans & \nu_1 -  {\bh_{11} }^\trans { H_{21} }^{-1} \bnu_2
        \end{array}
 \right].   
\end{equation*}
Hence, if $ \nu_1 - {\bh_{11} }^\trans { H_{21} }^{-1} \bnu_2 = 0$,
\begin{equation*}
  \rank \left[ \begin{array}{cc}
                      {\bh_{11} }^\trans & \nu_1  \\
                      H_{21} & \bnu_2
                  \end{array}
          \right] = j-1,
\end{equation*}
since $ \rank H_{21} = j-1$. Hence,
\[ \bnu = \left[ \begin{array}{c}
                         \nu_1 \\
                         \bnu_2
                       \end{array}
             \right] = \left[ \begin{array}{c}
                                    {\bh_{11} }^\trans \\
                                     H_{21}
                       \end{array}
             \right] \bs,
\]
where $\bs \ne \bze$. Then, 
\[ \bze = V_j^1 \bnu = V_j^1 \left[ \begin{array}{c}
                                    {\bh_{11} }^\trans \\
                                     H_{21}
                       \end{array}
             \right] \bs
          = V_j^1 H_j \left[ \begin{array}{c}
                                    {\rm I}_{j-1} \\
                                     \bze^\trans
                       \end{array}
             \right] \bs
           =A_{11} V_j^1 \left[ \begin{array}{c}
                                    {\rm I}_{j-1} \\
                                     \bze^\trans
                       \end{array}
             \right] \bs.
\]

Since $ A_{11}$ is nonsingular,
\[ V_j^1 \left[ \begin{array}{c}
                                    {\rm I}_{j-1} \\
                                     \bze^\trans
                       \end{array}
             \right] \bs = \left[ \bv_1^1,\ldots,\bv_{j-1}^1 \right] \bs = \bze,
\]
where $\bs \ne \bze$. Hence, $ \bv_1^1,\ldots,\bv_{j-1}^1 $ are linearly dependent and \\
$ \rank V_{j-1}^1 = \rank \left[ \bv_1^1,\ldots,\bv_{j-1}^1 \right] \le j-2$, but 
$ \rank V_j^1 = \rank \left[ \bv_1^1,\ldots,\bv_{j-1}^1, \bv_j^1 \right] = j-1 $. Hence, we have $ \rank V_{j-1}^1 = j-2 $.

Next, we will use an induction argument on $\ell$, where $ 1 \le \ell \le j-2$. 
Note 
\begin{equation}
    h_{i+1,i}\ne 0 \hspace{3mm} (1 \le i \le \ell).
 \label{hi+1i}
\end{equation}
Let $ \rank V_{\ell+1}^1 = \ell$ where $ V_{\ell+1}^1 \in {\bf R}^{r \times (\ell+1)} $.
Since $ \rank V_{\ell+1}^1 + \dim \nul (V_{\ell+1}^1 ) = \ell + 1$, we have 
$ \dim \nul (V_{\ell+1}^1) =1$.
Hence, let $ \nul ( V_{\ell+1}^1 ) = \left\{ c \bnu \right\} $, where $ c \in \real $, and
\[ \bnu = \left[ \begin{array}{c}
                          \nu_1 \\
                          \bnu_2
                       \end{array}
               \right] \ne \bze \in \real^{\ell+1}, \; \nu_1 \in \real, \; \bnu_2 \in {\bf R}^l.
\]
Noting that, $ A_{11} V_\ell^1 = V_{\ell+1}^1  \overline{H}_\ell $, similarly to (\ref{avvj1}), we have
\[ \br_\ell^1 = \beta \bv_1^1 - A_{11} V_\ell^1 \by = V_{\ell+1}^1 \left( \beta \be_1 - \overline{H}_\ell \by \right) , \]
where $ \be_1 = (1,0,\ldots,0)^\trans \in {\bf R}^{\ell+1} $.

Let   
\[ \overline{H}_\ell =  \left[ \begin{array}{c}
                                              {\bh_{11}}^\trans \\
                                              H_{21} 
                                  \end{array}
                          \right],
\]
where $ {\bh_{11}}^\trans = \left[ h_{11},\ldots, h_{1\ell} \right]$, and
\[ H_{21} = \left[ \begin{array}{ccc}
                            h_{21} & \cdots & h_{2\ell} \\
                                    & \ddots & \vdots \\
          $\mbox{\Large 0}$  &           & h_{\ell+1,\ell}
                          \end{array}
                  \right] ,
\]
where $ H_{21} $ is nonsingular due to (\ref{hi+1i}).

Then, note the following:
\begin{equation*}
 \begin{array}{lll}
 & \mbox{A least squares solution is obtained at step $\ell$ } \\
 \Longleftrightarrow & \exists \by \mbox{ such that } \br_\ell^1 = V_{\ell+1}^1 \left( \beta \be_1 - \overline{H}_\ell \by \right) = \bze \\
 \Longleftrightarrow & \exists \by \mbox{ such that } \beta \be_1 - \overline{H}_\ell \by \in \nul(V_{\ell+1}^1 ) \\
 \Longleftrightarrow & \exists \by \mbox{ such that } \left\{ \begin{array}{ccc}
                                                                          \beta -  {\bh_{11}}^\trans \by & = & c \nu_1 \\
                                                                          -H_{21} \by & = & c \bnu_2 
                                                                      \end{array} 
                                                              \right. \\
 \Longleftrightarrow & \nu_1 -  {\bh_{11}}^\trans { H_{21} }^{-1} \bnu_2 \ne 0 
 \end{array}
\end{equation*}

Hence, if $ \nu_1 -  {\bh_{11}}^\trans { H_{21} }^{-1} \bnu_2 \ne 0 $, a least squares solution is obtained at step $\ell$.

If $ \nu_1 -  {\bh_{11}}^\trans { H_{21} }^{-1} \bnu_2 = 0$, a least squares solution is not obtained at step $\ell$, and
\[ \left| \begin{array}{cc}
              \nu_1 &  {\bh_{11}}^\trans \\
              \bnu_2 & H_{21}
            \end{array}
    \right| =0 .
\]
Since $ H_{21} $ is nonsingular and $ \bnu \ne \bze $,
\[ \bnu = \left[ \begin{array}{c}
                         \nu_1 \\
                         \bnu_2
                      \end{array}
              \right] 
           =  \overline{H}_\ell \bs, 
\] 
where $ \bs \ne \bze \in {\bf R}^\ell $. Then,
\[ A_{11} V_\ell^1 \bs= V_{\ell+1}^1 \overline{H}_\ell \bs =  V_{\ell+1}^1  \bnu = \bze. \]
Since $ A_{11} $ is nonsingular, $ \rank V_\ell^1 \le \ell-1.$ But since $ \rank V_{\ell+1}^1 = \ell $, $\rank V_\ell^1 = \ell-1$.

Thus, by induction on $ \ell $, a least squares solution is obtained at step $ \ell \; ( 2 \le \ell \le j ) $, or $ \rank V_1^1 
= \rank \left[ \bv_1^1 \right] =0 $, so that $ \bv_1^1 = \bze$. Then, $ \br_1^1= \beta \bv_1^1 - A_{11} \bv_1^1 y = \bze$, 
so a least squares solution is obtained at step $1$.

Hence, if $ h_{i+1,i}\ne 0 \; (1 \le i \le j-1), \; h_{j+1,j}=0$, a least squares solution is obtained by step $j \; (j \le r+1)$.
\qed

The necessity of the condition $A_{12}=0$ for GMRES to determine a least squares solution of (\ref{lstsq}) for all $\bb, \bx_0 \in \rn$ was proved in Theorem 2.6 of Hayami and Sugihara\cite{HS}.


\begin{thebibliography}{10}
\bibitem{HS}
Hayami K, Sugihara M.
A geometric view of Krylov subspace methods on singular systems.
Numer Linear Algebra Appl. 2011; 18:449--469.

\bibitem{SS}
Saad Y, Schultz MH. 
GMRES: A generalized minimal residual algorithm for solving nonsymmetric linear systems. 
SIAM J Sci Statist Comput. 1986; 7:856--869.

\bibitem{BW}
Brown P, Walker HF. 
GMRES on (nearly) singular systems. 
SIAM J Matrix Anal Appl. 1997; 18:37--51.

\bibitem{SHZ}
Sugihara K, Hayami K, Zheng, N.
Right preconditioned MINRES for singular systems.
Numer Linear Algebra Appl. 2020; 27:e2277. https://doi.org/10.1002/nla.2277

\end{thebibliography}
\end{document}